\documentclass[12pt,a4paper]{amsart}

\usepackage{amsfonts}
\usepackage{amsmath}
\def\frak{\mathfrak}

\newtheorem{remark}{Remark}
\newtheorem{lemma}{Lemma}
\newtheorem{proposition}{Proposition}

\def\br{\mathbb R}
\def\bc{\mathbb C}

\def\bz{\mathbb Z}
\def\Ad{{\rm Ad}}
\def\ad{{\rm ad}}

\def\mand{\quad\mbox{and}\quad}

\begin{document}
\bibliographystyle{alpha}

\title{Stein extensions of Riemann symmetric spaces and some generalization}
\author[T. Matsuki]{Toshihiko MATSUKI\\
\bigskip {Dedicated to Professor Keisaku Kumahara on his 60th birthday}}
\address{\hskip-\parindent
        Toshihiko Matsuki\\
        Faculty of Integrated Human Studies\\
        Kyoto University\\
        Kyoto 606-8501, Japan}
\email{matsuki@math.h.kyoto-u.ac.jp}
\date{}

\begin{abstract}
It was proved by Huckleberry that the Akhiezer-Gindikin domain is included in the ``Iwasawa domain'' using complex analysis. But we can see that we need no complex analysis to prove it. In this paper, we generalize the notions of the Akhiezer-Gindikin domain and the Iwasawa domain for two associated symmetric subgroups in real Lie groups and prove the inclusion. Moreover, by the symmetry of two associated symmetric subgroups, we also give a direct proof of the known fact that the Akhiezer-Gindikin domain is included in all cycle spaces.
\end{abstract}

\maketitle

\section{Introduction}

\bigskip
\noindent 1.1 Akhiezer-Gindikin domain and Iwasawa domain

\bigskip
Let $G_\bc$ be a connected complex semisimple Lie group and $G_\br$ a connected real form of $G_\bc$. Let $K_\bc$ be the complexification in $G_\bc$ of a maximal compact subgroup $K$ of $G_\br$. Let $\frak{g}_\br=\frak{k}\oplus\frak{m}$ denote the Cartan decomposition of $\frak{g}_\br={\rm Lie}(G_\br)$ with respect to $K$. Let $\frak{t}$ be a maximal abelian subspace in $i\frak{m}$. Put
$$\frak{t}^+=\{Y\in\frak{t}\mid |\alpha(Y)|<{\pi\over 2} \mbox{ for all }\alpha\in\Sigma\}$$
where $\Sigma$ is the restricted root system of ${\frak g}_\bc$
with respect to $\frak{t}$. Then the ``Akhiezer-Gindikin domain'' $D$ in $G_\bc$ is defined in \cite{AG} by
$$D=G_\br(\exp\frak{t}^+)K_\bc.$$

Let $B$ be a Borel subgroup of $G_\bc$ such that $G_\br B$ is closed in $G_\bc$. Then $K_\bc B$ is the unique open dense $K_\bc$-$B$ double coset in $G_\bc$ (\cite{M3}). Define an open subset
$$\Omega=\{x\in G_\bc\mid x^{-1}G_\br B\subset K_\bc B\}$$
in $G_\bc$. Clearly, $\Omega$ is left $G_\br$-invariant and right $K_\bc$-invariant. The connected component $\Omega_0$ of $\Omega$ containing the identity is often called the ``Iwasawa domain''.

\begin{remark} \ {\rm Let $S_j\ (j\in J)$ be the $K_\bc$-$B$ double cosets in $G_\bc$ of (complex) codimension one and $T_j$ the closure of $S_j$. Then the complement of $K_\bc B$ in $G_\bc$ is
$\bigcup_{j\in J} T_j.$
So we can write
\begin{align*}
\Omega & =\{x\in G_\bc\mid x^{-1}G_\br\cap T_j=\phi\mbox{ for all }j\in J\} \\
& =\{x\in G_\bc\mid x\notin gT_j^{-1}\mbox{ for all }j\in J\mbox{ and }g\in G_\br\}.
\end{align*}
Thus $\Omega_0$ is Stein because $\Omega$ is the complement of an infinite family $\{gT_j^{-1}\mid j\in J,\ g\in G_\br\}$ of complex hypersurfaces.
}\end{remark}

The equality
\begin{equation}
D=\Omega_0 \tag{1.1}
\end{equation}
was proved in \cite{GM} when $G_\br$ is classical type or exceptional Hermitian type. Independently, Kr{\" o}tz and Stanton proved $D\subset \Omega_0$ for classical cases in \cite{KS}. We should note that the proofs in these two papers are based on elementary linear algebraic computations. (Remark that \cite{FH} did not refer \cite{GM} in their historical reference on (1.1).)
On the other hand, Barchini proved the inclusion
$D\supset \Omega_0$
by a general but elementary argument in \cite{B}.
Recently, Huckleberry (\cite{H}, Proposition 2.0.2 in \cite{FH}) gave a proof of the opposite inclusion
\begin{equation}
D\subset \Omega_0 \tag{1.2}
\end{equation}
by using strictly plurisubharmonicness of a function $\rho$ proved in \cite{BHH} (revised version).

But we can see that we need no complex analysis to prove (1.2). In this paper, we will generalize the notions of the Akhiezer-Gindikin domain and the Iwasawa domain for real Lie groups as in the next subsection and prove the inclusion (1.2).

\bigskip
\noindent 1.2 Generalization to real Lie groups

\bigskip
Let $G$ be a connected real semisimple Lie group and $\sigma$ an involution of $G$. Take a Cartan involution $\theta$ such that
$\sigma\theta=\theta\sigma.$
Put $H=G^\sigma_0=\{g\in G\mid \sigma(g)=g\}_0$ and $H'=G^{\sigma\theta}_0$. Then the symmetric subgroup $H'$ is called ``associated'' to the symmetric subgroup $H$ (c.f. \cite{M1}). The structure of the double coset decomposition $H'\backslash G/H$ is precisely studied in \cite{M4} in a general setting for arbitrary two involutions.

\begin{remark} \ {\rm (i) \ For example, if $G=G_\bc$, then $K_\bc$ and $G_\br$ are associated.

(ii) \ We should remark here that all the results on Jordan decompositions and elliptic elements etc. for the decomposition $G_\br\backslash G_\bc/K_\bc$ in Section 3 of \cite{FH} were already proved in \cite{M4}. In \cite{FH}, they are not referred at all.
}\end{remark}

Let $\frak{g}=\frak{k}\oplus\frak{m}=\frak{h}\oplus\frak{q}$ be the $+1,-1$-eigenspace decompositions of $\frak{g}={\rm Lie}(G)$ with respect to $\theta$ and $\sigma$, respectively. Then the Lie algebra $\frak{h}'$ of $H'$ is written as
$\frak{h}'=(\frak{k}\cap\frak{h})\oplus (\frak{m}\cap\frak{q}).$
Let $\frak{t}$ be a maximal abelian subspace of $\frak{k}\cap\frak{q}$. Then we can define the root space
$$\frak{g}_\bc(\frak{t},\alpha)=\{X\in\frak{g}_\bc\mid [Y,X]=\alpha(Y)X\mbox{ for all }Y\in\frak{t}\}$$
for any linear form $\alpha:\frak{t}\to i\br$. Here $\frak{g}_\bc=\frak{g}\oplus i\frak{g}$ is the complexification of $\frak{g}$. Put
$$\Sigma=\Sigma(\frak{g}_\bc,\frak{t})=\{\alpha\in i\frak{t}^*-\{0\} \mid \frak{g}_\bc(\frak{t},\alpha)\ne\{0\}\}.$$
Then $\Sigma$ satisfies the axiom of the root system (\cite{R}). Since $\theta(Y)=Y$ for all $Y\in\frak{t}$, we can decompose $\frak{g}_\bc(\frak{t},\alpha)$ into $+1,-1$-eigenspaces for $\theta$ as
\begin{equation}
\frak{g}_\bc(\frak{t},\alpha)=\frak{k}_\bc(\frak{t},\alpha) \oplus\frak{m}_\bc(\frak{t},\alpha). \tag{1.3}
\end{equation}
Define a subset
$\Sigma(\frak{m}_\bc,\frak{t})=\{\alpha\in i\frak{t}^*-\{0\} \mid \frak{m}_\bc(\frak{t},\alpha)\ne\{0\}\}$
of $\Sigma$ and put
$$\frak{t}^+=\{Y\in{\frak t}\mid |\alpha(Y)|<{\pi\over 2} \mbox{ for all }\alpha\in\Sigma(\frak{m}_\bc,\frak{t})\}.$$
Then we define a generalization of the Akhiezer-Gindikin domain $D$ in $G$ by
$$D=H'T^+H.$$
where $T^+=\exp\frak{t}^+$. (We will show in Proposition 1 that $D$ is open in $G$.)

Let $P$ be an arbitrary parabolic subgroup of $G$ such that $H'P$ is closed in $G$. Then $HP$ is open in $G$ by \cite{M3}. Define an open subset
$$\Omega=\{x\in G\mid x^{-1}H'P\subset HP\}=\{x\in G\mid H'x\subset PH\}$$
in $G$. Then we may call the connected component $\Omega_0$ a ``generalized Iwasawa domain''. As a complete generalization of (1.2), we can prove the following theorem.

\bigskip
\noindent {\bf Theorem.}\quad $D\subset \Omega_0$.

\bigskip
Here we should explain the construction of this paper to prove this theorem.
In Section 2, we prove properties on the generalized Akhiezer-Gindikin domain $D$. Lemma 2 is the most important basic technical lemma. It implies
$(\frak{h}'\cap \Ad(a)\frak{h})\oplus (\frak{q}'\cap \Ad(a)\frak{q})\subset\frak{k}$
for $a\in T^+$ where $\frak{q}'=(\frak{k}\cap\frak{q})\oplus (\frak{m}\cap\frak{h})$. Note that $H'\cap aHa^{-1}$ is the isotropy subgroup of the action of $H'$ at the point $aH\in G/H$. So the inclusion $\frak{h}'\cap \Ad(a)\frak{h}\subset\frak{k}$ is a generalization of Proposition 2 in \cite{AG}. But we do not find statements in \cite{AG} corresponding to the inclusion
$\frak{q}'\cap \Ad(a)\frak{q}\subset\frak{k}.$
This inclusion is the key of this paper.

In Proposition 1, we show that $D$ is open in $G$. This is a generalization of Proposition 4 in \cite{AG}. We also give a precise orbit structure
$H'\backslash D/H\cong T^+/W$
of $D$ where $W=N_{K\cap H}(\frak{t})/Z_{K\cap H}(\frak{t})$ in Proposition 2. This is a generalization of Proposition 8 in \cite{AG}.
In Section 3, we construct a left $H'$-invariant right $H$-invariant real analytic function $\rho$ on $D$ and prepare the key lemma (Lemma 3) which follows from Lemma 2.
In Section 4, we prove Theorem. Basic formulation is the same as Proposition 2.0.2 in \cite{FH}. But we do not need complex analysis.

\bigskip
\noindent 1.3 Application to cycle spaces

\bigskip
Note that we can exchange the roles of $H$ and $H'$. We can aplly this to the pair of $K_\bc$ and $G_\br$ in $G_\bc$. Let $P$ be a parabolic subgroup of $G_\bc$ such that $S=K_\bc P$ is closed in $G_\bc$. Then $S'=G_\br P$ is open in $G_\bc$. Put
$$\Omega(S)=\{x\in G_\bc\mid x^{-1}K_\bc P\subset G_\br P\}.$$
Then by Theorem, we have
$D^{-1}=K_\bc T^+G_\br \subset\Omega(S)_0$
by the notation in Section 1.1 and therefore
$D\subset\Omega(S)_0^{-1}.$
Since
$\Omega(S)^{-1} =\{x\in G_\bc\mid xK_\bc P\subset G_\br P\} =\{x\in G_\bc\mid xS\subset S'\},$
the domain $\Omega(S)_0^{-1}$ is usually called the ``cycle space'' for $S'$. Hence we have given a direct proof of the known fact:

\bigskip
\noindent {\bf Corollary.}\quad Akhizer-Gindikin domain $D$ is included in all cycle spaces.

\bigskip
\noindent (Remark: This fact was known by combining (1.2) and Proposition 8.3 in \cite{GM} because Proposition 8.3 implies that the Iwasawa domain $\Omega_0$ is included in all the cycle spaces $\Omega(S)_0^{-1}$.)

\bigskip
\noindent {\bf Aknowledgement}: The author would like to express his hearty thanks to S. Gindikin for his encouragement and for many useful discussions.

\section{A generalization of Akhiezer-Gindikin domain}

We will use the notations in Section 1.2 (not in Section 1.1). In this section, we will prepare basic results on a generalization of the Akhiezer-Gindikin domain by extending elementary arguments in \cite{M4} Section 3.

First we should note that we have only to consider the problem on each minimal $\sigma$-stable ideals of $\frak{g}$. So we may assume that $\frak{g}$ has no proper $\sigma$-stable ideals. We may also assume that $\frak{g}$ is noncompact, otherwise we have $P=G$ and the problem is trivial.

Let $[\frak{m},\frak{m}]$ denote the linear subspace of $\frak{k}$ spanned by $[Y,Z]$ for $Y,Z\in\frak{m}$. Then $\frak{g}'=[\frak{m},\frak{m}]\oplus\frak{m}$ becomes a $\sigma$-stable ideal of $\frak{g}$ because $\frak{g}'$ is $\sigma$-stable and
$[\frak{k},[\frak{m},\frak{m}]]\subset [[\frak{k},\frak{m}],\frak{m}] +[\frak{m},[\frak{k},\frak{m}]]\subset[\frak{m},\frak{m}].$
Since $\frak{g}'=\frak{g}$, we have $\frak{k}=[\frak{m},\frak{m}]$ and therefore
\begin{equation}
\frak{k}_\bc=[\frak{m}_\bc,\frak{m}_\bc]. \tag{2.1}
\end{equation}

\begin{lemma} \ {\rm (i)} \ Every root $\alpha\in\Sigma$ such that $\frak{k}_\bc (\frak{t},\alpha)\ne\{0\}$ is written as a sum of two elements in $\Sigma(\frak{m}_\bc,\frak{t})\cup\{0\}$.

{\rm (ii)} \ If $Y\in\frak{t}^+$, then $|\alpha(Y)|<\pi$ for all $\alpha\in\Sigma$.

{\rm (iii)} \ If $Y\in\frak{t}$ satisfies $\alpha(Y)=0$ for all $\alpha\in\Sigma$, then $Y=0$.

{\rm (iv)} \ If $Y,\ Z\in\frak{t}^+$ satisfy $\exp Y\in (\exp Z)(T\cap H)$, then $Y=Z$. Especially, $\exp$ is injective on $\frak{t}^+$.
\end{lemma}

Proof. \ (i) \ Let $X$ be a nonzero element in $\in\frak{k}_\bc (\frak{t},\alpha)$. Then by (2.1), we can write
$X=\sum_j [Y_j,Z_j]$
with some nonzero elements $Y_j\in\frak{m}_\bc (\frak{t},\beta_j)$ and $Z_j\in\frak{m}_\bc (\frak{t},\gamma_j)$ such that
$\beta_j+\gamma_j=\alpha.$

(ii) is clear by (i) and the definition of $\frak{t}^+$.

(iii) \ If $\alpha(Y)=0$ for all $\alpha\in\Sigma$, then $Y$ commutes with $\frak{g}_\bc (\frak{t},\alpha)$ for all $\alpha\in\Sigma\cup\{0\}$. Hence $Y$ commutes with $\frak{g}$ and therefore $Y=0$ because $\frak{g}$ is semisimple.

(iv) \ Suppose $\exp Y=(\exp Z)h$ with some $Y,\ Z\in\frak{t}^+$ and $h\in T\cap H$. Applying $\sigma$, we have
$\exp(-Y)=(\exp(-Z))h.$
So we have
$\exp 2Y=\exp 2Z.$
Put $a=\exp 2Y=\exp 2Z$ and apply $\Ad(a)$ to $\frak{g}_\bc (\frak{t},\alpha)$ for $\alpha\in\Sigma(\frak{m}_\bc,\frak{t})$. Then we have
$e^{2\alpha(Y)}=e^{2\alpha(Z)}.$
and therefore $\alpha(Y)-\alpha(Z)\in \pi i\bz$. Since $|\alpha(Y)|<\pi/2$ and $|\alpha(Z)|<\pi/2$ by the definition of $\frak{t}^+$, we have
$\alpha(Y-Z)=\alpha(Y)-\alpha(Z)=0\mbox{ for all }\alpha\in\Sigma(\frak{m}_\bc,\frak{t}).$
This implies $\alpha(Y-Z)=0\mbox{ for all }\alpha\in\Sigma$ by (i) and therefore $Y-Z=0$ by (iii).
\hfill q.e.d.

\bigskip
Let $Y$ be an element of $\frak{t}$ and put $a=\exp Y$. Consider the conjugate
$\sigma_a=\Ad(a)\sigma\Ad(a)^{-1}=\sigma\Ad(a)^{-2}$
of $\sigma$ by $\Ad(a)$. Then $\sigma_a$ is an involution of $\frak{g}$ such that
$$\frak{g}^{\sigma_a}=\Ad(a)\frak{h}\mand G_0^{\sigma_a}=aHa^{-1}.$$
Put $\tau=\sigma\theta=\theta\sigma$. The key idea of \cite{M4} was to consider the automorphism $\tau\sigma_a$ (which is not involutive in general) of $\frak{g}$. Since
$\tau\sigma_a=\tau\sigma\Ad(a)^{-2}=\theta\Ad(a)^{-2}$
and since $\theta$ and $\Ad(a)$ commutes, the automorphism $\tau\sigma_a$ is semisimple. Hence by Lemma 1 in \cite{M4}, we have a direct sum decomposition
\begin{equation}
\frak{g}=(\frak{h}'+\Ad(a)\frak{h})\oplus (\frak{q}'\cap\Ad(a)\frak{q}) \tag{2.2}
\end{equation}
where $\frak{q}'=(\frak{k}\cap\frak{q})\oplus (\frak{m}\cap\frak{h})$. On the other hand, we have the $+1,-1$-eigenspace decomposition of $\frak{g}^{\tau\sigma_a}$ for $\tau$ by
$$\frak{g}^{\tau\sigma_a}=(\frak{g}^{\tau\sigma_a}\cap \frak{h}')\oplus (\frak{g}^{\tau\sigma_a}\cap \frak{q}')
=(\frak{h}'\cap \Ad(a)\frak{h})\oplus (\frak{q}'\cap \Ad(a)\frak{q})$$
as in \cite{M4} Section 3.

\begin{lemma} \ If $Y\in\frak{t}^+$, then $\frak{g}^{\tau\sigma_a}=\frak{z}_\frak{k}(Y)=\{X\in\frak{k}\mid [Y,X]=0\}$.
\end{lemma}

Proof. \ By (1.3), we have a direct sum decomposition
\begin{equation}
\frak{g}_\bc=\bigoplus_{\alpha\in\Sigma\cup\{0\}} (\frak{k}_\bc(\frak{t},\alpha)\oplus \frak{m}_\bc(\frak{t},\alpha)). \tag{2.3}
\end{equation}
If $X\in \frak{k}_\bc(\frak{t},\alpha)$, then we have
$\tau\sigma_a X=\theta\Ad(a)^{-2} X=e^{-2\alpha(Y)}X.$
On the other hand, if $X\in \frak{m}_\bc(\frak{t},\alpha)$, then we have
$\tau\sigma_a X=\theta\Ad(a)^{-2} X=-e^{-2\alpha(Y)}X.$
Hence the decomposition (2.3) is the eigenspace decomposition of $\frak{g}_\bc$ for $\tau\sigma_a$. We have only to verify whether every direct summand in (2.3) is contained in $\frak{g}_\bc ^{\tau\sigma_a}$ or not.

Note that $\alpha(Y)$ is pure imaginary. Let $X$ be a nonzero element in $\frak{m}_\bc(\frak{t},\alpha)$. Since $|\alpha(Y)|<\pi/2$ by the assumption, we have
$-e^{-2\alpha(Y)}\ne 1.$
Hence $X\notin \frak{g}_\bc ^{\tau\sigma_a}$. On the other hand, let $X$ be a nonzero element in $\frak{k}_\bc(\frak{t},\alpha)$. Since $|\alpha(Y)|<\pi$ by Lemma 1 (ii), we have
$e^{-2\alpha(Y)}=1\Longleftrightarrow \alpha(Y)=0.$
Hence
$X\in \frak{g}_\bc ^{\tau\sigma_a}\Longleftrightarrow \alpha(Y)=0.$
Thus we have proved
$$\frak{g}_\bc ^{\tau\sigma_a}=\bigoplus_{\alpha\in\Sigma\cup\{0\},\ \alpha(Y)=0} \frak{k}_\bc(\frak{t},\alpha)=\frak{z}_{\frak{k}_\bc}(Y)$$
and therefore $\frak{g}^{\tau\sigma_a}=\frak{z}_\frak{k}(Y)$. \hfill q.e.d.

\bigskip
Since $\frak{g}^{\tau\sigma_a}=\frak{z}_\frak{k}(Y)$ is a compact Lie algebra and since $\frak{t}$ is maximal abelian in $\frak{q}'\cap \Ad(a)\frak{q}=\frak{z}_{\frak{k}\cap\frak{q}} (Y)$, we have
$\frak{q}'\cap \Ad(a)\frak{q}=\Ad(H'\cap aHa^{-1})_0 \frak{t}.$
Moreover if $U$ is a neighborhood of the origin $0$ in $\frak{t}$, then
\begin{equation}
V=\Ad(H'\cap aHa^{-1})_0 U\mbox{ is a neighborhood of $0$ in }\frak{q}'\cap \Ad(a)\frak{q}. \tag{2.4}
\end{equation}

\begin{proposition} \ $D$ is open in $G$.
\end{proposition}

Proof. \ By the left $H'$-action and the right $H$-action on $D$, we have only to show that a neighborhood of $a=\exp Y$ for $Y\in\frak{t}^+$ is contained in $D$. Take a neighborhood
$U=\{Z-Y\mid Z\in\frak{t}^+\}$
of $0$ in $\frak{t}$. Then we have $T^+a^{-1}=\exp U$. Hence
\begin{align}
e & \in H'T^+Ha^{-1}=H'T^+a^{-1}aHa^{-1}=H'(\exp U)aHa^{-1}. \tag{2.5} \\
& =H'\exp(\Ad(H'\cap aHa^{-1})_0U)aHa^{-1}=H'(\exp V)aHa^{-1}. \notag
\end{align}
Since $V$ is a neighborhood of $0$ in $\frak{q}'\cap \Ad(a)\frak{q}$ by (2.4), it follows from (2.2) that $H'T^+Ha^{-1}$ contains a neighborhood of $e$. Hence $H'T^+H$ contains a neighborhood of $a$. \hfill q.e.d.

\begin{proposition} \ Let $a$ and $b$ be elements of $T^+$. Then
$b=\ell ah^{-1}\mbox{ for some }\ell\in H'\mbox{ and }h\in H\Longleftrightarrow b=waw^{-1}\mbox{ for some }w\in N_{K\cap H}(\frak{t}).$
Here $N_{K\cap H}(\frak{t})$ is the normalizer of $\frak{t}$ in $K\cap H$. $($Write $K=G^\theta$ as usual.$)$
\end{proposition}

Proof. \ Since the implication $\Longleftarrow$ is clear, we have only to prove $\Longrightarrow$. Suppose $b=\ell ah^{-1}$ for some $\ell\in H'$ and $h\in H$. Put
$\frak{t}'=\Ad(\ell)\frak{t}=\Ad(\ell a)\frak{t}=\Ad(bh)\frak{t}.$
Then $\frak{t}'$ is a maximal abelian subspace of $\frak{q}'\cap\Ad(b)\frak{q}$. Since $(\frak{h}'\cap\Ad(b)\frak{h},\ \frak{q}'\cap\Ad(b)\frak{q})$ is a compact symmetric pair by Lemma 2, there is an $x\in H'\cap bHb^{-1}$ such that
$\Ad(x)\frak{t}'=\frak{t}.$
Put $\ell'=x\ell\in H'$ and $h'=b^{-1}xbh\in H$. Then we have
$$\ell'a{h'}^{-1}=x\ell ah^{-1}b^{-1}x^{-1}b=xbhh^{-1}b^{-1}x^{-1}b=b$$
and
\begin{align}
\ell'T{h'}^{-1} & =x\ell Th^{-1}b^{-1}x^{-1}b
=x\ell aTh^{-1}b^{-1}x^{-1}b \tag{2.6} \\
& =xbhTh^{-1}b^{-1}x^{-1}b
=xT'x^{-1}b
=Tb=T \notag
\end{align}
where $T=\exp\frak{t}$ and $T'=\exp\frak{t}'$.

Since $\ell'{h'}^{-1}=\ell'e{h'}^{-1}\in T$ by (2.6) and since
$\ell'T{\ell'}^{-1}\ell'{h'}^{-1}=T$
also by (2.6), we have
$\ell'T{\ell'}^{-1}=T.$
Write $\ell'=k\exp X$ with $k\in K\cap H$ and $X\in\frak{m}\cap\frak{q}$. Then it is well-known that $[X,\frak{t}]=\{0\}$ by a standard argument. (Suppose $\Ad(k\exp X)Y=Y'$ for some $Y,\ Y'\in\frak{k}$. Then applying $\theta$, we have $\Ad(k\exp(-X))Y=Y'$ and therefore $\Ad(\exp 2X)Y=Y$. Since $\ad X: \frak{g}\to\frak{g}$ is expressed by a real symmetric matrix, we have $[X,Y]=0$.)

Hence we have $kTk^{-1}=T$ and
$b=\ell'a{\ell'}^{-1}\ell'{h'}^{-1}
=kak^{-1}\ell'{h'}^{-1}.$
We have only to prove $c=\ell'{h'}^{-1}=e$. Write $h'=k'\exp X'$ with $k'\in K\cap H$ and $X'\in\frak{m}\cap\frak{h}$. Then it follows from $\ell'=ch'$ that
$k\exp X=ck'\exp X'.$
Since $c\in T\subset K$, we have $k=ck'$ and $X=X'=0$. Hence
$c\in T\cap H.$
Since $kak^{-1},\ b\in\exp{t}^+$, we have $c=e$ by Lemma 1 (iv). \hfill q.e.d.

\section{Construction of a function $\rho$}

Write $W=N_{K\cap H}(\frak{t})/Z_{K\cap H}(\frak{t})$. Let $\rho_0$ be a $W$-invariant real analytic function on $\frak{t}^+$ which has no critical points except the origin. For the sake of later use, we should also assume $\rho_0(Z)$ tends to $+\infty$ when $Z$ goes to the boundary of $\frak{t}^+$. For example, we may put
\begin{equation}
\rho_0(Z)=\sum_{\alpha\in\Sigma(\frak{m}_\bc,\frak{t})} {1\over \pi-2i\alpha(Z).} \tag{3.1}
\end{equation}
This function (3.1) is clearly convex and therefore it has no critical points except the origin. By Proposition 2 and Lemma 1 (iv), we can define a function $\rho$ on $D$ by
$$\rho(\ell(\exp Z)h)=\rho_0(Z)\quad\mbox{ for }\ell\in H',\ h\in H\mbox{ and }Z\in\frak{t}^+.$$

\begin{proposition} \ $\rho$ is real analytic on $D$.
\end{proposition}

Proof. \ By the left $H'$-action and the right $H$-action on $D$, we have only to show $\rho$ is real analytic at every $a=\exp Y\in T^+$. Consider the right $a$-translate
$$\rho_a(x)=\rho(xa)\quad \mbox{for }x\in Da^{-1}$$
of $\rho$. Since $\rho_a$ is left $H'$-invariant and right $aHa^{-1}$-invariant, we have only to show that the function
$\rho'_a(X)=\rho_a(\exp X)\quad\mbox{for }X\in V$
is real analytic at $0$ by (2.2) and (2.5) where $U\subset\frak{t}$ and $V=\Ad(H'\cap aHa^{-1})_0U\subset\frak{q}'\cap \Ad(a)\frak{q}$ are as in the proof of Proposition 1. Since 
$\rho_a(gxg^{-1})=\rho_a(x)$
for $g\in H'\cap aHa^{-1}$ and $x\in \exp V$, $\rho'_a$ is $\Ad(H'\cap aHa^{-1})$-invariant. Note that $(H'\cap aHa^{-1})_0=(Z_{K\cap H}(Y))_0$ by Lemma 2 and that
$$\rho'_a(Z)=\rho((\exp Z)a)=\rho_0(Z+Y)\quad\mbox{for }Z\in U$$
is invariant under the action of $w\in W$ such that $w(Y)=Y$. Since we can easily extend the well-known Chevalley's restriction theorem to real analytic functions at $0$, the function $\rho'_a$ on $V=\Ad(H'\cap aHa^{-1})_0U$ is real analytic at $0$. \hfill q.e.d.

\begin{remark} \ {\rm (i) \ The function $\rho$ on $D$ has no critical points outside $H'H=H'eH$ by the assumption on $\rho_0$.

(ii) \ If $\rho_0$ is a $W$-invariant smooth function on $\frak{t}^+$. Then we can show that $\rho$ is a smooth function on $D$ by using \cite{S}.
}\end{remark}

The tangent space $T_a(G)$ of $G$ at $a$ is identified with $\frak{g}=T_e(G)$ by the right $a$-action. In other words, we identify $T_a(G)$ with the left infinitesimal action of $\frak{g}$ at $a$. Now we have the following key lemma.

\begin{lemma} \ Let $a=\exp Y$ with $Y\in\frak{t}^+-\{0\}$. Then the hyperplane in $T_a(G)$ defined by $d\rho=0$ is orthogonal, with respect to the Killing form on $\frak{g}$, to a nonzero vector $Z$ in $\frak{k}$.
\end{lemma}

Proof. \ Taking the right $a$-translate $\rho_a$ of $\rho$ as in the proof of Proposition 3, we have only to consider the hyperplane in the tangent space $T_e(G)\cong \frak{g}$ defined by $d\rho_a=0$.

Since $\rho_a$ is left $H'$-invariant and right $aHa^{-1}$-invariant, the differential $d\rho_a$ vanishes on $\frak{h}'+\Ad(a)\frak{h}$. Hence the normal vector $Z$ is contained in the orthogonal complement $\frak{q}'\cap \Ad(a)\frak{q}\ (\subset\frak{k}$ by Lemma 2) of $\frak{h}'+\Ad(a)\frak{h}$. \hfill q.e.d.

\section{Proof of Theorem}

Proof of Theorem. \ Basic formulation is the same as Proposition 2.0.2 in \cite{FH}. Suppose that
$D\not\subset PH.$
We will deduce a contradiction.

Let $PxH$ be a $P$-$H$ double coset with the least dimension among the $P$-$H$ double cosets intersecting $D$. So the intersection $PxH\cap D$ is relatively closed in $D$. Since $H'P=(K\cap H)P$ by \cite{M3}, we have $H'=(K\cap H)(P\cap H')=(P\cap H')(K\cap H)$. Hence $PxH$ intersects $(K\cap H)T^+$ and the image of $\rho|_{PxH\cap D}$ is equal to the image of $\rho|_{PxH\cap (K\cap H)T^+}$. The set $\{x\in (K\cap H)T^+\mid \rho(x)\le m\}$ is compact for any $m\in\br$ because we carefully assumed that $\rho_0$ is $+\infty$ on the boundary of $\frak{t}^+$. Hence the function $\rho|_{PxH\cap D}$ attains its minimum on some point $ka$ with $k\in K\cap H$ and $a\in T^+$. Replacing $P$ by the $k$-conjugate $k^{-1}Pk$, we may assume $k=e$. Since $a\in PxH$ and $PxH\cap PH=\phi$, we have $a\ne e$.

By Lemma 3, there is a nonzero element $Z$ in $\frak{k}$ such that $Z$ is orthogonal to $\frak{p}={\rm Lie}(P)$. But this leads a contradiction because $Z\in\frak{k}$ is also orthogonal to $\theta\frak{p}$ and therefore
$Z\mbox{ is orthogonal to }\frak{p}+\theta\frak{p}=\frak{g}$
which cannot happen since the Killing form is nondegenerate on $\frak{g}$. \hfill q.e.d.

\end{document}